\documentclass[12pt]{article}
\begin{document}

{\Large\bf Conditional Linearizability Criteria for Scalar Fourth
Order Semi-Linear Ordinary Differential Equations}

F M Mahomed$^1$ and A Qadir$^2$

$^1$School of Computational and Applied Mathematics, Centre for
Differential Equations, Continuum Mechanics and Applications\\
University of the Witwatersrand, Wits 2050, South Africa\\
Fazal.Mahomed@wits.ac.za\\
$^2$Centre for Advanced Mathematics \& Physics\\
National University of Sciences \& Technology\\
Campus of College of Electrical \& Mechanical Engineering\\
Peshawar Road, Rawalpindi, Pakistan\\
aqadirmath@yahoo.com\\

{\bf Abstract}. Using geometric methods for linearizing systems of
second order cubically semi-linear ordinary differential equations
and third order quintically semi-linear ordinary differential
equations, we extend to the fourth order by differentiating the
third order equation. This yields criteria for linearizability of a
class of fourth order semi-linear ordinary differential equations,
which have not been discussed in the literature previously. It is
shown that the procedure can be extended to higher order. Though the
results for the higher orders are complicated, they are doable by
algebraic computing. The standard Lie approach, as developed at
present does not seem to be amenable to giving results that can be
handled even by algebraic computing.


\section{Introduction}
First order ordinary differential equations (ODEs) can always be
linearized (i.e. converted to linear form) [1] by point
transformations [2]. Lie [3] showed that all second order ODEs that
can be converted to linear form must be cubically semi-linear
\begin{equation}
y''+E_3(x,y)y'^3+E_2(x,y)y'^2+E_1(x,y)y'+E_0(x,y)=0, \label{(1)}
\end{equation}
the coefficients $E_0$ to $E_3$ satisfying the over-determined integrable system
\begin{eqnarray}
b_x&=&-\frac13 E_{1y}+\frac23 E_{2x}+be-E_0E_3,\nonumber\\
b_y&=&E_{3x}-b^2+bE_2-E_1E_3+eE_3,\nonumber\\
e_x&=&E_{0y}+e^2-eE_1-bE_0+E_0E_2,\nonumber\\
e_y&=&\frac23E_{1y}-\frac13E_{2x}-be+E_0E_3,\label{(2)}
\end{eqnarray}
where $b$ and $e$ were (for Lie's purposes) auxiliary variables and
the suffices $x$ and $y$ refer to partial derivatives. The
auxilliary variables arise naturally in the geometric approach
mentioned shortly.

The linearization programme was carried forward to the third order
by Chern [4, 5] using contact transformations to discuss the
linearizability of equations reducible to the linear forms
$u'''(t)=0$ and $u'''(t)+u(t)=0$ and Grebot [6, 7] used point
transformations that were restricted to the class of transformations
$t=\phi(x)$, $u=\psi(t,x)$ for the same purpose. It was subsequently
shown [8] that there are three classes of third order ODES that are
linearizable by point transformations, viz. those that reduce to the
above two forms or $u'''(t)+ \alpha (t) u(t)=0$. (Note that one
cannot choose $\alpha$ to be a constant, as it would then reduce to
the second form by re-scaling the independent variable.) After that
Neut and Petitot [9] and later Ibragimov and Maleshko [10] used the
original Lie procedure [3], of postulating some point transformation
that yields the desired linearizability, and determining the
criteria for it to exist. The resulting criteria are very
complicated and barely manageable for algebraic computing. Clearly,
higher order attempts would not be feasible for this method.
Maleshko [11] also provided a simple algorithm to reduce third order
odes of the form $y'''=F[y,y',y'']$ to second order and used the
$\partial/\partial t$ symmetry to reduce the order to two. He then
used the Lie linearizability criteria to determine the
linearizability of the third order equation. This could be trivially
extended to the fourth order for $y''''=F[y',y'',y''']$, but this
extension would not work for any non-trivial problem.

Though no progress was made towards providing explicit procedures
for linearizing systems of ODEs, it was shown [12] that there will
generally be multiple classes, even for second order systems. The
number of classes for even dimensions, $2m$, is $2m^2+3$ and for odd
dimensions, $(2m-1)$ it is $2m^2-2m+4$. Clearly, there was no
question of proceeding to higher order systems with the previous
procedures.

Using the connection between the isometry algebra and the symmetries
of the system of geodesic equations [13, 14], linearizability
criteria were stated for a system of second order quadratically
semi-linear ODEs, of a class that could be regarded as a system of
geodesic equations [15], that we will call {\it of geodesic type}.
The criteria came from requiring that the coefficients in the
equations, regarded as Christoffel symbols, yield a zero curvature
tensor. The criteria also provide the required compatibility
conditions. Further, even for larger dimensional systems, it is
possible to associate a metric with the system of geodesics when
these criteria are met [16]. The flatness of the metric allows a
coordinate transformation to be defined from the given metric to a
Cartesian (Eulidean or pseudo-Euclidean) form. These are the
linearizing transformations. They can be determined using complex
variables [15] and used to write down the solution in terms of the
original variables.

Utilizing the projection procedure of Aminova and Aminov [13], which
appeals to the fact that the geodesic equations do not depend on the
geodetic parameter, the system of $n$ second order ODEs of geodesic
type can be reduced to system of  $(n-1)$ second order cubically
semi-linear ODEs [17]. When applied to a system of two dimensions we
get a scalar cubically semi-linear ode that yields the Lie cirteria!
Applied to a system of three dimensions, one obtains a system of two
cubically semi-linear odes with extended Lie cirteria. The system of
odes so obtained is in the class of maximally symmetric equations
(out of the total of five classes mentioned earlier).

Differentiating the quadratically and cubically semi-linear system
of ODEs relative to the independent variable gives third order ODEs.
Taking the general class of the scalar third order ODE one gets
linearizability criteria for scalar third order ODEs [18]. {\it This
class is is not included in the Neut and Petitot} [9] {\it and
Ibragimov and Maleshko classes} [10]. Though there can be an overlap
with the Maleshko class [11] it is obviously not contained in that
either. It may be wondered how it can be outside the three classes
allowed by [8]. The reason is that it is {\it not} obtained by point
transformations but is linearized {\it conditionally} to a second
order ODE being linearizable. In effect its first integral is a
linearizable second order semi-linear ODE. This is, thus, a new type
of linearizability.

The quadratically semi-linear system of geodesic type, on
differentiation yields a cubically semi-linear third order system of
the same dimension. This is obtained by replacing the second order
term in the differentiated equation by a quadratic expression using
the original equation. It could also be written as a third order ODE
with a second order term multiplied by a first order term and then
followed by a term quadratic in the first derivative. The
linearizability criteria for this system for two dimensions have
also been derived [19].

It is clear that the geometrical approach has far-reaching
consequences for the linearizability of ODEs. In this paper it will
be used for stating criteria for the conditional linearizability of
fourth order ODEs in some detail aand its utilization for higher
orders will also be discussed. First a brief summary of the notation
used, relevant for our present purposes, will be provided.

\section{Notation and review}
Equations of geodesic type for two dimensions are of the form
\begin{eqnarray}
x''=a(x,y)x'^2+2b(x,y)x'y'+c(x,y)y'^2, \nonumber\\
y''=d(x,y)x'^2+2e(x,y)x'y'+f(x,y)y'^2. \label{(3)}
\end{eqnarray}
There are six coefficients and they may be identified with the six
Christoffel symbols for a 2-dimensional space,
\begin{eqnarray}
\Gamma_{11}^1=-a, \Gamma^1_{12}=-b, \Gamma_{22}^1=-c, \nonumber\\
\Gamma_{11}^2=-d, \Gamma^2_{12}=-e, \Gamma_{22}^2=-f. \label{(4)}
\end{eqnarray}
Requiring that the space corresponding to these Christoffel symbols
is flat
\begin{equation}
R^{i}_{jkl}=\Gamma^i_{jl,k}-\Gamma^i_{jk,l}+\Gamma^i_{mk}\Gamma^m_{jl}-\Gamma^i_{ml}
\Gamma^m_{jk}=0, \label{(5)}
\end{equation}
yields the linearizability conditions
\begin{eqnarray}
a_y-b_x + be -cd = 0 , b_y -c_x+ (ac - b^2) + (bf - ce) = 0 , \nonumber\\
d_y-e_x - (ae - bd) - (df - e^2) = 0 , (b + f)_x = (a + e)_y ,
\label{(6a)}
\end{eqnarray}
which provide the metric coefficients through
\begin{eqnarray}
p_x = -2(ap + dq) , q_x = -bp - (a + e)q - dr , r_x = -2(bq + er) , \nonumber\\
p_y = -2(bp + eq) , q_y = -cp - (b + f)q - er , r_y = -2(cq + fr) . \label{(7)}
\end{eqnarray}
The compatibility of these equations is guaranteed by the curvature
tensor being zero. The metric coefficients being obtained the
linearizing transformation is available as a coordinate
transformation and its inverse can then be computed.

The above system of two equations can be projected to the scalar
cubically semi-linear ODE
\begin{equation}
y''+cy'^3-(f-2b)y'^2+(a-2e)y'-d=0, \label{(8)}
\end{equation}
which is clearly of the form (1) for suitable values of the
coefficients in the former equation. It is here that one sees the
``auxilliary variables" arising naturally. Note that there are only
four coefficients here as against the six Christoffel symbols and
for the original system of two equations. This degeneracy allows for
extra classes when this procedure is applied to larger systems of
ODEs but does not interfere with the uniqueness for the scalar
equation [16]. We shall write $g=f-2b$ and $h=a-2e$ for convenience.
Again, the coordinate transformations can be obtained, but there is
freedom of choice of $g$ and $h$ that can lead to more or less
convenient coordinate transformations. The linearizability criteria
can then be given in the form [20]
\begin{eqnarray}
3(ch)_x+3dc_y-2gg_x-gh_y-3c_{xx}-2g_{xy}-h_{yy}=0 , \nonumber\\
3(dg)_y+3cd_x-2hh_y-hg_x-3d_{yy}-2h_{xy}-g_{xx}=0. \label{(9)}
\end{eqnarray}
We shall be using this form.

This equation can be differentiated to yield a third order equation
of the general form
\begin{equation}
y'''+(A_2y'^2-A_1y'+A_0)y''+B_4y'^4-B_3y'^3+B_2y'^2-B_1y'+B_0=0,
\label{(10)}
\end{equation}
which is linear in the second derivative, but with a coefficient
quadratic in the first derivative, where the coefficients are of
this equation are identified with the original equation by,
\begin{equation}
c=A_2/3, g=A_1/2, h=A_0, \label{(11)}
\end{equation}
\begin{equation}
d=-\int B_2dx+k(y)=\int (B_1-A_{0x})dy+l(x), \label{(12)}
\end{equation}
the constant of integration being fixed by compatibility with the
metric tensor, provided the compatibility conditions
\begin{equation}
B_4=A_{2y}/3, B_3=A_{1y}/2-A_{1x}/3, B_2=A_{0y}-A_{1x}/2,
\label{(13)}
\end{equation}
hold and the linearizability conditions (6) are satisfied. Once
again, the coordinate transformations are obtainable and provide the
solution.

The above third order ODE is a total derivative. Though it may not
be immediately obvious that a given ODE is an exact derivative, the
equation in this form may seem trivial in some sense. One can use
(8) in (10) to replace the second derivative term and obtain a third
order ODE quintically semi-linear in the first derivative,
\begin{equation}
y'''-\alpha y'^5+\beta y'^4-\gamma y'^3+\delta y'^2-\epsilon
y'+\phi=0, \label{(14)}
\end{equation}
where
\begin{equation}
\alpha=3c^2, \beta=5cg+c_y, \gamma=4ch+2g^2+g_y-c_x,
\delta=3cd+3gh+h_y-g_x, \label{(15)}
\end{equation}
and the compatibility conditions are
\begin{equation}
\epsilon=2dg+h^2+d_y-h_x, \phi=dh-d_x. \label{(16)}
\end{equation}
{\it This} equation is {\it not} an exact derivative, and hence
provides a non-trivial utilization of the procedure given. The
coefficients can be inverted to obtain the original coefficients and
hence the metric coefficients, yielding the solution.

We shall use the above formulation to obtain the linearizability
criteria for higher order equations in the next section.

\section{Conditional linearizability criteria for fourth order ODEs}
On differentiating (14) for the scalar third order ODE, writing the
independent variable as $x$and the dependent variable as $y$, we get
\begin{eqnarray}
y''''-(5\alpha y'^4-4\beta y'^3+3\gamma y'^2-2\delta y'+\epsilon )y''
- \alpha y'^6 +(\beta _y-\alpha _x)y'^5-(\gamma _y-\beta _x)y'^4+ \nonumber\\
(\delta _y-\gamma _x)y'^3-(\epsilon _y-\delta _x)y'^2 + (\phi_
x-\epsilon _x)y'+\phi _x=0. \label{(17)}
\end{eqnarray}
The general form of this equation is
\begin{equation}
y''''-(A_4y'^4-A_3y'^3+A_2y'^2-A_1y'+A_0)y''-B_6y'^6+B_5y'^5-B_4y'^4
+B_3y'^3-B_2y'^2+B_1y'-B_0=0 , \label{(18)}
\end{equation}
subject to the identification of coefficients
\begin{eqnarray}
c = \sqrt{A_4/15}, g = (A_3-4c_y)/20c, h = (A_2-6g^2-3g_y+3c_x)/12c,
\nonumber\\ d = (A_1-6gh-2h_y+2g_x)/6c, \label{(19)}
\end{eqnarray}
with the constraints
\begin{eqnarray}
B_6=A_4/5, B_5=A_{3y}/4-A_{4x}/5, B_4=A_{2y}/3-A_{3x}/4,
B_3=A_{1y}/2-A_{2x}/3, \nonumber\\
B_2=A_{0y}-A_{1x}/2, B_1=dh_y+hd_y-d_{xy}-A_{0x}/2,
B_0=dh_x+hd_x-d_{xx}, \nonumber\\ A_0=2gd+h^2+d_y-h_x. \label{(20)}
\end{eqnarray}
In this form it is a total derivative and may be regarded as
trivial in some sense.

Replacing the second derivative by the first derivative expressions
from (8) we get the fourth order ODE semi-linear in the first
derivative in the seventh degree
\begin{equation}
y''''+P_7y'^7-P_6y'^6+P_5y'^5-P_4y'^4+P_3y'^3-P_2y'^2+P_1y'-P_0=0,
\label{(21)}
\end{equation}
which is {\it not} a total derivative. For consistency of the
identification of coefficients
\begin{eqnarray}
c = (P_7/15)^{1/3}, g = P_6/35c^2-2c_y/7c, \nonumber\\
h = P_5/27c^2-26g^2/27c+c_x/3c-gc_y/3c^2-8g_y/27c-c_{yy}/27c^2,
\nonumber\\
d = P_4/21c^2-38gh/21c-2g^3/7c^2+8gc_x/21c^2-8hc_y/21c^2 \nonumber\\
+g_x/3c-gg_y/3c^2-2h_y/7c+2c_{xy}/21c^2-g_{yy}/21c^2, \label{(22)}
\end{eqnarray}
so that all four coefficients are explicitly given. We now have
four differential constraints
\begin{eqnarray}
P_3=28cdg+13ch^2+12g^2h-3(h+d)c_x+4dc_y-(2g+3h)g_x+(3h+2d)g_y \nonumber\\
-(c+3g)h_x+2(g+h)h_y-3cd_x+(c+2g)d_y+g_{xx}-2h_{xy}+d_{yy}, \nonumber\\
P_2=18chd+8g^2d+7gh^2-6dc_x-5hg_x+4dg_y-4gh_x+4hh_y-3cd_x+2gd_y \nonumber\\
+g_{xx}-2h_{xy}+d_{yy}, \nonumber\\
P_1=6cd^2-8ghd+h^3-4dg_x-3hh_x+3dh_y-2gd_x+2hd_y+h_{xx}-2d_{xy}, \nonumber\\
P_0=2gd^2+h^2d-2dh_x-hd_x+dd_y+d_{xx}, \label{(23)}
\end{eqnarray}
apart from the earlier stated linearizability conditions (9).

Instead of starting with the quintically semi-linear third order
ODE, we can start with the form involving the second derivative
(10). In this case we get
\begin{eqnarray}
y''''+(A_2y'^2-A_1y'+A_0)y'''+(B_1y'-B_0)y''^2+(C_3y'^3-C_2y'^2+C_1y'-C_0)y'' \nonumber\\
-(D_5y'^5-D_4y'^4+D_3y'^3-D_2y'^2+D_1y'-D_0)=0 , \label{(24)}
\end{eqnarray}
which is semilinear involving the third third derivative with a
coefficient quadratic in the first derivative, second derivative
squared with a coefficient linear in the first derivative, and
otherwise quintic in the first derivative. It is a total derivative
with the identification
\begin{equation}
c = A_2/3=B_1/6, g = A_1/2=B_0/2, h = A_0, \label{(25)}
\end{equation}
the differential constraints
\begin{eqnarray}
C_3 = 7A_{2y}/3, C_2 = 5A_{1y}/2-2A_{2x}, C_1 = 3A_{0y}-2A_{1x},
\nonumber\\
D_5=A_{2yy}/3, D_4=A_{1yy}/2-2A_{2xy}/3,
D_3=A_{0yy}/2-A_{1xy}+A_{2xx}/3, \label{(26)}
\end{eqnarray}
and the additional constraints that yield $d$:
\begin{equation}
C_0=d_y-2A_{0x}, D_0=d_{xx}, \label{(27)}
\end{equation}
which gives the value as a single and a double integral to be solved
simultaneously
\begin{equation}
d= \int (C_0+2A_{0x})dy+k(x)= \int(\int D_0 dx)dx+l(y)x+m(y),
\label{(28)}
\end{equation}
where $k, l, m$ are arbitrary functions of the respective variables.
We have the additional requirements for $D_1$ and $D_2$,
\begin{equation}
D_1+2C_{0x}+3A_{0xx}=0, D_2-C_{0y}-A_{1xx}/2=0. \label{(29)}
\end{equation}

One can now use the third order equation (10) to replace the third
derivative term in (24) and obtain the equation in a form
quadratically semi-linear in the second derivative that has a
coefficient quadratically semilinear in the first derivative for the
higher order, quartic in the first derivative for the one linear in
the second derivative and otherwise of sixth order in the first
derivative,
\begin{eqnarray}
y''''+(Q_1y'-Q_0)y''^2-(R_4y'^4-R_3y'^3+R_2y'^2-R_1y'+R_0)y''
\nonumber\\
-(S_6y'^6-S_5y'^5+S_4y'^4-S_3y'^3+S_2y'^2-S_1y'+S_0)=0 ,
\label{(30)}
\end{eqnarray}
with the identification of the first three coefficients given by
\begin{equation}
c=Q_1/6, g=Q_0/2, h=(R_2-Q_0^2+Q_{1x}-5Q_{0y}/2)/Q_1. \label{(31)}
\end{equation}
Since the coefficient $h$ is very complicated we shall use it as a
shorthand for the above expression. The fourth coefficient is given
by the solution of
\begin{equation}
d=\int (R_0-h^2+2h_x)dy+k(x). \label{(32)}
\end{equation}
In the subsequent equations replacing $d_y$ by the integrand above
and writing $d_x$ (as there is no convenient expression for it), the
following constraint equations must hold
\begin{eqnarray}
R_4=Q_1^2/4, R_3=Q_1Q_0+7Q_{1y}/6, R_1=2Q_0h+h^2+3h_y-2Q_{0x}, \nonumber\\
S_6=Q_1Q_{1y}/12, S_5=-Q_1Q_{1x}/36+Q_0Q_{1y}/6+Q_1Q_{0y}/4+Q_{1yy}/6, \nonumber\\
S_4=-Q_0Q_{1x}/6+hQ_{1y}/6-Q_1Q_{0x}/4+Q_0Q_{0y}/2+Q_1h_y/2-Q_{1xy}/3+Q_{0yy}/2,
\nonumber\\
S_3=Q_1(R_0-h^2-h_x)/2-hQ_1/6-Q_0Q_{0x}/2+hQ_{0y}/2+Q_0h_y+Q_{1xx}/6-Q_{0xy}+h_{yy},
\nonumber\\
S_2=Q_1d_x/2+Q_0(R_0-h^2+h_x)-hh_y-hQ_{0x}/2+R_{0y}+Q_{0xx}/2, \nonumber\\
S_1=h(R_0-h^2+5h_x)-2R_{0x}-Q_0d_x-3h_{xx}, S_0=hd_x-d_{xx},
\label{(33)}
\end{eqnarray}
along with the linearization criteria. This is not a total
derivative. Though looking more messy, it is perfectly usable.

Yet another form that is not a total derivative can be obtained from
(24) by using (8) to replace the second derivative term. In this
case we get
\begin{eqnarray}
y''''+(A_2y'^2-A_1y'+A_0)y'''+B_7y'^7-B_6y'^6+B_5y'^5-B_4y'^4+B_3y'^3 \nonumber\\
-B_2y'^2+B_1y'-B_0=0 , \label{(34)}
\end{eqnarray}
with the identification of all four coefficients given by
\begin{eqnarray}
c=\sqrt{A_2/3}, g=A_1/2, h=A_0,  \nonumber\\
d=(B_4-8cA_1A_0-\frac{1}{4}A_1^3+3A_1c_x-7A_0c_y+2cA_{1x} \nonumber\\
-5A_1A_{1y}/4-3cA_{0x}+2c_{xy}-\frac{1}{2}A_{1yy})/4A_2,\label{(35)}
\end{eqnarray}
where we have largely used the symbol $c$ rather than $A_2$ to avoid
the extra complications due to a square root and its
differentiation. Further, though $d$ is given algebraically, it is
too complicated to use conveniently. The constraint equations are
\begin{eqnarray}
B_7=6c^3,B_6=7c(A_1+c_y),B_5=3A_2A_0+5cA_1^2/2-6cc_x+7A_1c_y/2+5cA_{1y}/2+c_{yy},
\nonumber\\
B_3=8cA_1d+6cA_0^2+A_1^2A_0-6A_0c_x+7dc_y-A_1A_{1x}+5A_0A_{1y}/2
\nonumber\\
-2cA_{0x}+3A_1A_{0y}/2+cd_y+c_{xx}-A_{1xy}+A_{0yy}, \nonumber\\
B_2=12cA_0d+A_1^2d+A_1A_0^2-6dc_x-2A_0A_{1x}+5dA_{1y}/2-A_1A_{0x})
\nonumber\\
+3A_0A_{0y}+A_1d_y/2+A_{1xx}/2-2A_{0xy}+d_{yy}, \nonumber\\
B_1=6cd^2+2A_1A_0d-4dA_{1x}-2A_0A_{0x}+3dA_{0y}+A_0d_y+A_{0xx}-2d_{xy},
\nonumber\\
B_0=d(A_1d+d_y-A_{0x})+d_{xx}, \label{(36)}
\end{eqnarray}
along with the linearization conditions. Again, this seems very
messy but it has the advantage that all four coefficients are
identified.

We could have one other procedure, to use both equations (8) and
(10) to replace the second and third derivatives. Since this only
involves the first derivative to the same power as before, this will
not yield a new class. Thus we have the following theorems.

{\bf Theorem 1:} {\it Equation (21) is linearizable with the
identifications (22) if the constraints (23) and the linearizability
criteria (9) are satisfied, where ($c \ne 0$), with $h=a-2e,
g=f-2b$, after requiring consistency of the Christoffel symbols with
the deduced metric coefficients.}

It is to be noted that the four linearizability conditions (6) are
stated in terms of the 6 coefficients $a,...,f$ and not the four
coefficients $c,d,g,h$. Thus there is degeneracy in the choices
available. Any choices of $a$ and $e$ for a given $h$, or $f$and $b$
for a given $g$, compatible with the metric coefficient relations
(7) are permissible. For each such choice we would get corresponding
linearizability conditions. Instead, we use the conditions (9) to
check linearizability and exploit the freedom of choice to construct
a convenient metric.

{\bf Theorem 2:} {\it Equation (30) is linearizable with the
identifications (31) and (32) if the constraints (33) along with the
linearizability criteria (9) are satisfied, after requiring
consistency of the Christoffel symbols with the deduced metric
coefficients.}

{\bf Theorem 3:} {\it Equation (34) is linearizable with the
identifications (35) if the constraints (36)  and the
linearizability criteria (9) are satisfied, after requiring
consistency of the Christoffel symbols with the deduced metric
coefficients.}

Note that there are two other equations, (18) and (24), that are
total derivatives of linearizable equations with appropriate
identifications, that are also linearizable subject to the
corresponding constraints and conditions. We do not stress on these
because they may be regarded, in some sense, as trivial.

\section{Examples}
In this section we present some examples of fourth order equations
that can be linearized by our procedure.

{\bf 1.} The equation
\begin{equation}
y''''-(18y'^2/y^2-16y'/y+k^2-5l)y''+12y'^4/y^3 -8ky'^3/y^2+kly'=0 ,
\label{(37)}
\end{equation}
is a fourth order total derivative of the form of (18). However, it
is clearly not trivial to spot this fact by looking at the equation.
It satisfies the required constraints (20) and linearizability
conditions (9).

{\bf 2.} The equation
\begin{equation}
y''''-24y'^4/y^3+33ky'^3/y^2+(28l-10k^2)y'^2/y+(k^2-35l)ky'/2-(k^2-5l)ly=0,
\label{(38)}
\end{equation}
which is a fourth order equation that is not a total derivative, of
the form (21) and is linearizable as it satisfies the constraints
(23) and the linearizability conditions (9).

{\bf 3.} The equation
\begin{equation}
y''''-(4y'/y-k)y'''-4y''^2/y+(10y'^2/y^2-l)y'' -4y'^4/y^3=0 ,
\label{(39)}
\end{equation}
which is a total derivative of the form of (24) and may be more
easily identified as a total derivative. It is linearizable as it
satisfies the constraints (26) and (29) and the linearizability
conditions (9).

{\bf 4.} The equation
\begin{equation}
y''''-4y''^2/y^4-(6y'^2/y^2-8ky'/y+k^2+l)y''+4y'^4/y^3-2ky'^3/y^2-4kly'^2/y+kly'=0
, \label{(40)}
\end{equation}
which is of the form (30) and satisfies the constraints (33) and the
linearizability conditions (9). It is not a total derivative.

{\bf 5.} The equation
\begin{equation}
y''''-(4y'y-k)y'''+4y'^4/y^3+3ky'^3/y^2+(4l-k^2)y'^2/y-7kly'/2-3l^2y=0
, \label{(41)}
\end{equation}
which is not a total derivative and is of the form (34) satisfying
the constraints (36) and the linearizability criteria (9). As such
it is linearizable.

Though it is not apparent from looking at the equations they come
from differentiating the linearizable third order equation
\begin{equation}
y'''-6y'^3/y^2+8ky'^2/y-(k^2-5l)y'+kly=0, \label{(42)}
\end{equation}
and the total derivative linearizable equation
\begin{equation}
y'''-(4y'/y-k)y''+2y'^3/y^2-ly'=0, \label{(43)}
\end{equation}
(in some cases subject to conditions of some equation holding) which
both come from differentiating the second order linearizable
equation
\begin{equation}
y''-2y'^2/y+ky'/2+ly=0, \label{(44)}
\end{equation}
and in one case applying a condition. As such, their solutions
necessarily have two arbitrary constants but may not have more. It
may be noted that these equations do not have any explicit
dependence on $x$. That symmetry being guaranteed only one more
needs to be looked for to find the solution by symmetry analysis.
Using the procedure of writing the equation as of geodesic type, we
can write down the solution directly.

The fact that all of these fourth order equations have a common root
becomes obvious from our analysis right in the beginning, as the
identification of the four coefficients in each case gives
\begin{equation}
 c=0, g=2, h=k/2, d=-l. \label{(45)}
\end{equation}

{\bf 6.} The fourth order ODE
\begin{equation}
 y''''+ (3xy'^2+2/x)y'''+6xy'y''^2+ (6y'^2-4/x^2)y''+4y'/x^3=0, \label{(46)}
\end{equation}
is of the form of (24) and is hence a total derivative (though this
fact is not obvious by inspection), with the coefficients given by
(25), (27) and (28). It satisfies the constraints (26) and (29) and
the linearizability criteria (9). As such it is linearizable.

{\bf 7.} The fourth order ODE
\begin{equation}
 y''''+6xy'y''^2-(9x^2y'^4+6y'^2+8/x^2)y''+4y'/x^3=0, \label{(47)}
\end{equation}
is of the form of (30), and hence is not a total derivative, with
the coefficients given by (31) and (32) and satisfies the
constraints (33) and the linearizability criteria (9). As such it is
linearizable.

{\bf 8.} The fourth order ODE
\begin{equation}
 y''''+ (3xy'^2+2/x)y'''+6x^3y'^7+ 18y'^5+16y'^3/x+12y'/x^3=0, \label{(48)}
\end{equation}
is of the form of (34), and hence is not a total derivative, with
the coefficients given by (35) and satisfying the constraints (36)
and the linearizability criteria (9). As such it is linearizable.

{\bf 9.} The fourth order ODE
\begin{equation}
 y''''-(15x^2y'^4+21y'^2+6/x^2)y''-6xy'^5+12y'/x^3=0, \label{(49)}
\end{equation}
is of the form of (18), and hence is a total derivative, with the
coefficients given by (19), satisfying the constraints (20) and the
linearizability criteria (9). As such it is linearizable.

{\bf 10.} The fourth order ODE
\begin{equation}
 y''''+ (3xy'^2+2/x)y'''+6x^3y'^7+ 18y'^5+16y'^3/x+12y'/x^3=0, \label{(50)}
\end{equation}
is of the form of (21), and hence is not a total derivative, with
the coefficients given by (22) and satisfying the constraints (23)
and the linearizability criteria (9). As such it is linearizable.

It is again apparent that even the total derivative equation is not
obviously so, and that the other equations being linearizable would
not be clear by inspection. The coefficients of examples (6) to (10)
are all $c=x, h=2/x, g=d=0$ and the {\it root equation} is the
geodesic equation for flat space in polar coordinates. Thus the
linearizing transformation is simply the conversion from polar to
Cartesian coordinates.

It will be noticed that in the above examples either $x$ or $y$ are
missing from the coefficients. This makes the application of the
identification and constraints relatively trivial. We end with a
couple of non-trivial examples, in which both arise.

{\bf 11.} The fourth order ODE
\begin{equation}
 y''''-15x^3y'^7/y^6-15xy'^6/y^5+39y'^5/y^4+39y'^4/xy^3-36y'^3/xy^2-36y'^2/x^2y
 +24y'/x^3=0, \label{(51)}
\end{equation}
is of the form of (21), and hence is not a total derivative. Its
coefficients are given by (22) and come out to be $c=-x/y^2, g=1/y,
h=2/x, d=0$. They satisfy the constraints (23) and linearizability
criteria (9). As such this equation is linearizable. In fact it can
be derived by differentiating the linearizable third order equation
[19]
\begin{equation}
y'''-3x^2y'^5/y^4-3xy'^4/y^3+6y'^3/y^2+6y'^2/xy-6y'/x^2=0,
\label{(52)}
\end{equation}
writing the root equation from the coefficients, and using it to
substitute the $y''$ arising from the differentiation. The solution
is
\begin{equation}
Axy+Bx/y=1, \label{(53)}
\end{equation}
where $A$ and $B$ are arbitrary constant real numbers.

{\bf 12.} The fourth order ODE
\begin{eqnarray}
 y''''-(6xy'/y^2+2/y)y''^2-(9x^2y'^4/y^4-2xy'^3/y^3-7y'^2/y^2-8y'/xy+8/x^2)y''
 \nonumber \\
 +6x^2y'^6/y^5-2xy'^5/y^4-2y'^4/y^3-6y'^3/xy^2-4y'^2/x^2y+8y'/x^3=0, \label{(54)}
\end{eqnarray}
is of the form of (30). Its coefficients are given by (31) and (32)
and are the same as of the above example. Hence they share a common
root equation. This equation arises by differentiating the root
equation twice and using its first derivative to replace the $y'''$.
It satisfies the constraints (33) and the linearizability criteria
(9) and is linearizable, yielding the same linear equation and
possessing the same solution as of the previous example.

\section{Concluding remarks}

We have provided a procedure to determine the linearizability of
some fourth order scalar semi-linear odes that are linearizable. We
have written them as five classes. There could have been other ways
of getting to the five classes, e.g. by first differentiating and
then replacing or first replacing and then differentiating. Since,
the general form of both procedures would be the same they are not
different. Again, we could have used a lower order equation to {\it
partially} replace terms in the higher order equation. For example,
we could have used the second order equation to replace the
quadratic term in the second derivative in (52) but retain the
linear term as it is. These are not {\it really independent} in some
sense. In this sense there are five independent classes. Two of
these may be regarded as trivial as they are total derivatives.
However, a glance at the equations will show that even these are not
so easy to identify as total derivatives. The other three classes
are {\it not} total derivatives. All can be thought of as arising
from the same second order linearizable differential equation that
satisfies the Lie conditions, by double differentiation, but
conditional to the original differential equation and in some cases
the differentiated differential equation. As such, the
linearizability classes are {\it non-classical} and would not lie in
the three classical classes reducible to the fourth order odes
$y''(x)=0$, $y''(x)+y(x)=0$ and $y''+\alpha (x)y(x)=0$. These are
guaranteed four arbitrary constants while our solutions are only
guaranteed {\it two}.

Let us call the underlying second order equation the {\it root
equation}, and the specific fourth order equations {\it forms} of
the {\it similar} fourth order equations. It is clear that all forms
arising from a common second order equation will form an equivalence
class. As such, all similar fourth order equations will have two
solutions in common but may have other, different solutions. The
question arises whether these equivalence classes of linearizable
fourth order semi-linear odes are disjoint and can be used to
decompose the space of those odes that are linearizable of this
type. In this context, it seems unlikely that this be the case
because there would be some of them that are linearizable and have
two solutions in common with others. As such, we conjecture that the
set of linearizable fourth order odes with common root equations
will not be decomposable into disjoint classes.

It is worth noting that there is only one class of linearizable
second order semi-linear ode (as proved by Lie) and two classes of
third order odes, one of which is a total derivative. For the fourth
order we find {\it five} of which two are total derivatives and
three are not. Can this procedure be carried further? The answer is
``obviously it can". The fifth order will have one form that only
involves the first derivatives apart from the fifth derivative and
is of ninth order in them. It is obvious that here all four
coefficients will be easily identified, but now instead of the four
constraints (22) we will have {\it six} constraints. Similarly, for
the sixth order the corresponding class will have the first
derivative to the eleventh power and will have eight constraints and
so on. Of course the usual linearizability conditions (6) would also
have to be met.

The question arises as to the number of classes that we can now
have. Clearly, for the fifth order there will be five classes of
total derivative forms obtained by differentiating the five classes
of the fourth order. How many other classes will there be? The
number can be obtained by counting the different forms that one
could obtain and turns out to be six. Thus the total number of
classes is eleven. It is obvious that there will then be eleven
classes of total derivative sixth order equations obtained. One
could continue the procedure of counting but it would be nice to
have a general formula for the total number of classes for any
order.

It is worth stressing that as can be seen from the examples, our
procedure provides the solutions of the linearizable semi-linear
equations of our classes for any order that is identified. It should
not be too difficult to prepare an algebraic code for identification
of coefficients, the constraints and the linearizability conditions
for the various classes. Certainly, the form involving only the
first order could be relatively trivially provided. It would be
worth while to get the algebraic codes prepared for at least some of
the higher order classes.

{\bf Acknowledgements}\\
AQ is most grateful to DECMA and the School of Computational and
Applied Mathematics, University of the Witwatersrand.

\section*{References}
1. Some text book that gives this fact.\\
2. Lie S.  {Theorie der Transformationsgruppen}. {\it Math. Ann.}
1880, {\bf 16}: 441. \\
3. Lie S. Klassifikation und Integration von gew\"onlichen
Differentialgleichungenzwischen $x$, $y$, die eine Gruppe von
Transformationen gestaten. {\it  Arch. Math.} 1883, {\bf VIII, IX}: 187. \\
4. Chern, S.S., `Sur la geometrie d'une equation differentielle du
troiseme orde', {\it C.R. Acad. Sci. Paris}, (1937) 1227. \\
5. Chern, S.S., `The geometry of the differential equation
$y'''=F(x,y,y,y'')$', Sci, Rep. Nat. Tsing Hua Univ. 4 (1940), 97-111. \\
6. Grebot, G., `The linearization of third order ODEs', preprint 1996. \\
7. Grebot, G., `The characterization of third order ordinary
differential equations admitting a transitive fibre-preserving point
symmetry group', {\it J. Math. Anal. Applic.} {\bf 206}, (1997), 364-388. \\
8. Mahomed F M and Leach P G L.Symmetry Lie Algebras of $n$th Order
Ordinary Differential Equations. {\it J. Math Anal Applic} 1990, {\bf151}: 80. \\
9.Neut, S., and Petitot M., `La g\'eom\'etrie de l'\'equation
$y'''=f(x,y,y',y'')$', {\it C.R. Acad. Sci. Paris S\'er I} {\bf 335}, (2002), 515-518.\\
10. Ibragimov, N.H., and Meleshko, S.V., `Linearization of
third-order ordinary differential equations by point and contact
transformations' {\it J. Math. Anal. Applic.} {\bf 308}, (2005), 266-289.\\
11. Meleshko, S.V., `On linearization of third-order ordinary differential
equations' {\it J. Phys. A.: Math. Gen. Math.} {\bf 39}, (2006), 15135-45.\\
12. Wafo Soh, C., and Mahomed, F.M., `Symmetry breaking for a system
of two linear second-order ordinary differential equations',
{\it Nonlinear Dynamics}  {\bf22}, (2000), 121. \\
13. Aminova, A. V. and  Aminov, N. A.-M., `Projective geometry of
systems of differential equations: general conceptions',
{\it Tensor N S} {\bf 62}, 2000, 65-86.\\
14. Feroze, T., Mahomed, F. M. and Qadir, A., `The connection between
isometries and symmetries of geodesic equations of the underlying spaces',
{\it Nonlinear Dynamics} {\bf 45}, 2006, 65. \\
15. Mahomed, F. M., and Qadir, A., `Linearization criteria for a system
of second-order quadratically semi-linear ordinary differential equations',
{\it Nonlinear Dynamics}, {\bf 48}, 2007, 417.\\
16. Fredericks, E., Mahomed, F.M., Momoniat, E. and Qadir, A.,
`Constructing the space from a system of geodesic equations',
preprint, University of the Witwatersrand, Johannesburg, 2007. \\
17. Mahomed, F. M., and Qadir, A., `Linearization Criteria for Systems of
Cubically Semi-Linear Second-Order Ordinary Differential Equations',
preprint, University of the Witwatersrand, Johannesburg, 2007. \\
18. Mahomed, F.M. and Qadir, A., Linearizability criteria for a
class of third order semi-linear ordinary differential equations,
preprint 2007. (Presented in the Workshop in Honour of Leonhard
Euler, Centre for Differential Equations, Continuum Mechanics and
Applications, School of Computational and Applied Mathematics,
University of the Witwatersrand, Johannesberg, South Africa August 3, 4, 2007.) \\
19. Mahomed, F.M., Naeem, I. and Qadir, A., Linearizability of
semilinear systems of third order ordinary differential equations,
paper under preparation. \\
20. Tresse, A., `Sur les Invariants Diff\'erentiels des Groupes
Continus de Transformations', {\it Acta Math.} {\bf 18}, (1894), 1.

\end{document}